\documentclass{amsart}

\usepackage{epsfig}
\usepackage{graphicx}
\usepackage{color}

\usepackage{amssymb}
\usepackage{amsthm}
\usepackage{amsfonts}
\usepackage{amsmath,amssymb,amscd}
\usepackage{lineno,hyperref}

\usepackage{enumerate}

\newcommand{\CV}{\mathcal{C}(V)}
\newcommand{\CVp}{\mathcal{C}(V')}
\newcommand{\Lq}{\mathcal{L}_{q}}

\newcommand{\go}{\mathcal{G}_{\lambda,\omega}}

\renewcommand{\H}{\mathcal{H}}
\renewcommand{\P}{\mathcal{P}}
\newcommand{\F}{\mathcal{F}}
\newcommand{\G}{\mathcal{G}}
\renewcommand{\L}{\mathcal{L}}
\newcommand{\K}{\mathcal{K}}

\renewcommand{\AA}{\sf{A}}
\newcommand{\BB}{\sf{B}}
\newcommand{\DD}{\sf{D}}

\newcommand{\Ss}{\sf{S}}
\newtheorem{lemma}{Lemma}
\newtheorem{prop}{Proposition}
\newtheorem{theorem}{Theorem}
\newtheorem{corollary}{Corollary}

\def\qed{\hfill {\hbox{\footnotesize{$\Box$}}}}

\newcommand{\reals}{\mathbb{R}}

\def\BBox{\kern  -0.2cm\hbox{\vrule width 0.15cm height 0.3cm}}

\oddsidemargin .cm \evensidemargin 1.cm \textwidth=16.5cm
\include {mak}
\begin{document}

\oddsidemargin 16.5mm
\evensidemargin 16.5mm

\thispagestyle{plain}
\hspace{-.5cm}
\includegraphics[scale=0.5]{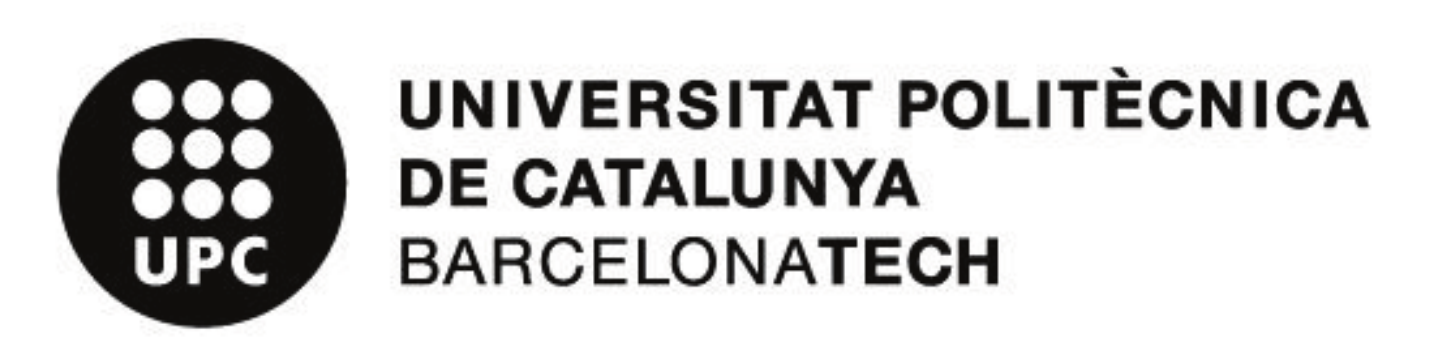} \hspace{.5cm}
\includegraphics[scale=0.14]{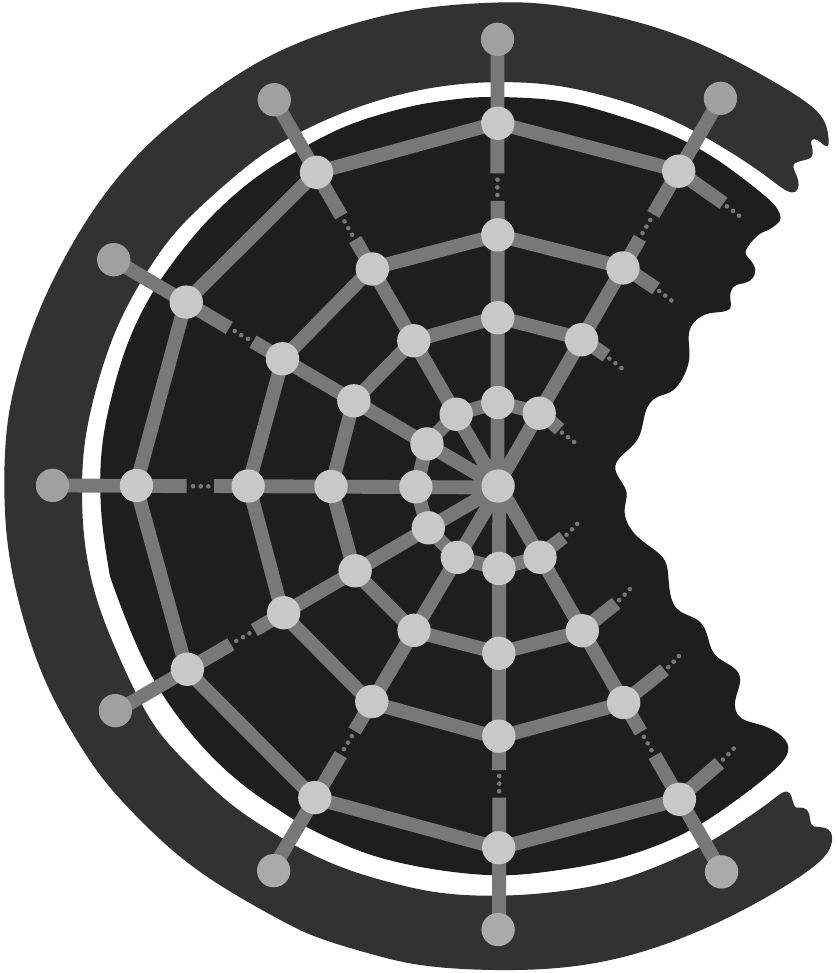}
\vspace{-1.3cm}

\hfill\begin{minipage}{.28\textwidth}
\noindent {\small\bf GRUPO MAPTHE}


\noindent {\scriptsize \bf Preprint 15.}

\noindent {\scriptsize \bf MAY 2015.}
\end{minipage}

\vspace{5cc}
\begin{center}

{\Large\bf  Green Operators of Networks with a new vertex
\rule{0mm}{6mm}\renewcommand{\thefootnote}{}
\footnotetext{\scriptsize 2000 Mathematics Subject Classification:  }}

\vspace{1cc} {\large\it A. Carmona, A.M. Encinas, S. Gago, M. Mitjana }

\vspace{1cc}
\parbox{24cc}{{\small{\bf Abstract.}
Any elliptic operator defines an automorphism on the orthogonal subspace to the eigenfunctions associated with the lowest eigenvalue, whose inverse is the orthogonal Green operator. In this study, we show that elliptic Schr\"{o}dinger operators on networks that have been obtained by adding a new vertex to a given network, can be seen as perturbations of the Schr\"{o}dinger operators on the initial network. Therefore, the Green function on the new network can be computed in terms of the Green function of the original network.
}}
\end{center}

\vspace{1cc}


\section{Introduction}
\label{introducton}

Discrete elliptic operators can be seen as the discrete counter part of elliptic partial differential operators. In particular, positive semi--definite Schr\"{o}dinger operators defined on a finite network are examples of those self--adjoint operators. Any elliptic operator defines an automorphism on the orthogonal subspace to the eigenfunctions associated with the lowest eigenvalue, whose inverse is the orthogonal Green operator. 

In \cite{CEM14}, some of the authors analyzed the effect of a
perturbation of the network by computing the effective resistance of the perturbed
networks through Sherman--Morrison--Woodbury like--formulas, instead of using
the Sherman--Morrison formula recursively. In fact, since adding edges
to a network does not modify the space of functions on the vertex set of the
network, this class of perturbation was placed into the general framework of
perturbations of discrete elliptic operators. Specifically, we showed that this
problem corresponds with the superposition of rank one perturbations that
are orthogonal to the eigenfunction associated with the lowest eigenvalue of
the elliptic operator.

The scenario changes when the perturbation consists on adding new vertices
to the network. Only few works have tackled the problem of adding a new vertex, see for instance\cite{B85}. In this work, we consider perturbations
that consist on adding a new vertex to a network. After some well--known
operations on the Schr\"odinger operator of the perturbed network, that involves the inverse of the Schur complement of the block corresponding to the added vertices, we show that this Schur complement can be seen as a perturbation of the Schr\"odinger operator of the original network, understood as a discrete elliptic operator, that is a superposition of rank one perturbations that, this time, are not orthogonal to the eigenfunction associated with the lowest eigenvalue of the elliptic
operator. Therefore, we can apply the general theory developed in \cite{CEM14} for this kind of perturbations.

We start the study by revisiting the perturbation of an elliptic operator with a sum of projections that can be, or not, orthogonal to the eigenfunction associated with the smallest eigenvalue.\ Thus, we consider the relation between the Green o\-pe\-rator of the new operator in terms of the Green operator of the previous one.

Next section is devoted to the application of the mentioned results to the addition of a new vertex to the network $\Gamma$ in order to get a network $\Gamma'$. Moreover, we obtain the relation between the   Schr\"{o}dinger operators of the two networks $\Gamma$ and $\Gamma'$ and in addition, we give the explicit expression of the matrix associated with the Green operator.

\section{Specific notation and preliminary results}
\label{}

Given a finite set $V$ of $n$ elements, we denote by $\CV$ the space of real valued functions on $V$. For any function $u\in \CV$, the associated vector in $\reals^n$ will be denoted by ${\sf u}$. For any vertex $x \in V$, the Dirac function at $x$ is denoted by $\varepsilon_x\in \CV$; the scalar product on $\CV$ is $\langle u,v \rangle = \sum_{x\in V} u_xv_x$ for each $u, v\in \CV$. A unitary and positive function $\omega$ is called a weight and $\Omega(V)$ denote the set of weights.

If $\K$ is an endomorphism of $\CV$, it is self-adjoint when $\langle \K(u),v\rangle=\langle u,\K(v)\rangle$ for any $u\in \CV$. Moreover, $\K$ is positive semi-definite when $\langle \K(u),u\rangle \geq 0$ for any $u\in \CV$. A self-adjoint operator $\K$ is elliptic if it is positive semi-definite and its lowest eigenvalue $\lambda$ is simple. Moreover, there exists a unique unitary function $\omega \in \CV$, up to sign, satisfying $\K(\omega)=\lambda \omega$, so $\K$ is called $(\lambda,\omega)$-elliptic operator. It is straighforward that a $(\lambda,\omega)$-elliptic operator is singular iff $\lambda=0$.  We denote by $\lambda^\dag=\lambda^{-1}$ iff $\lambda\neq 0$, and $\lambda^\dag=0$ otherwise.

Any function $K:V\times V\rightarrow \mathbb{R}$ is called a kernel on $V$, and it determines an endomorphim of  $\CV$ by assigning to any $u\in \CV$ the function $\K(u)=\sum_{y\in V}K(\cdot,y)u(y)$. Conversely, each endomorphism of $\CV$ is determined by the kernel $K(x,y)=\langle \K(\epsilon_y),\epsilon_x\rangle$ for any $x$, $y\in V$. Therefore, an endomorphism $\K$ is self-adjoint iff its kernel $K$ is a symmetric function.

Given $\sigma$, $\tau\in \CV$, we denote by $\P_{\sigma,\tau}$ the endomorphism of $\CV$ that assigns to each $u\in \CV$ the function $\P_{\sigma,\tau}(u)=\langle \tau, u\rangle \sigma$, and it is called a projector, as it assigns to any $u\in \CV$ its projection on $\sigma$ along $\tau$. Observe that the corresponding kernel is $P_{\sigma,\tau}(x,y)=(\sigma \otimes \tau)(x,y)=\sigma(x)\tau(y)$. In particular, when $\omega\neq 0$ the endomorfism $\P_{\omega,\omega}$ is denoted simply by $P_{\omega}$.

Given $\lambda\geq0$, $\omega\in\Omega(V)$ and a $(\lambda,\omega)$-elliptic operator $\F$, we will be concerned with the so-called Poisson equation for $\F$ on $V$: for a given $f\in \CV$ find $u\in \CV$ such that $\F(u)=f$. As $\F$ defines an automorphism on $\omega^\perp$, the inverse of a  $(\lambda,\omega)$-elliptic operator $\F$ on $\omega^{\perp}$ is called orthogonal Green operator and it is denoted by $\G$. This operator on $\omega^{\perp}$ can be extended to $\CV$ by assigning to any $f\in \CV$ the unique solution of the Poisson equation $\F(u)=f-\P_{\omega}(f)$.

Now consider a non-null function $\sigma\in \CV$, the associated self-adjoint projection $\P_{\sigma}$ and the operator $\H_{\sigma}=\F+\P_{\sigma}$, called perturbation of $\F$ by $\sigma$.
The relation between Green operators of both $\H_\sigma$ and $\F$ can be found  in  \cite[Corollary 3.6]{CEM14}.




\begin{corollary} \label{corol1}
Consider $\sigma\in \CV$. If $\lambda=\langle \sigma,  \omega \rangle=0$, then
$$\G^{\H_\sigma}= \G - \frac{1}{1+\langle \G(\sigma),\sigma \rangle}\P_{\G(\sigma)},$$
 whereas when    either $\lambda>0$ or $\langle \sigma,  \omega \rangle\neq 0$, then $\H_\sigma$ is invertible and
\begin{eqnarray*}
\H^{-1}_\sigma&=&\G- \frac{1}{\beta}
\Big[\lambda \P_{\G(\sigma)} - \langle \sigma,\omega \rangle \left(\P_{\G(\sigma),\omega} -  \P_{\omega,\G(\sigma)}\right)  -\big(1+\langle \G(\sigma),\sigma\rangle\big)\P_{\omega} \Big],
\end{eqnarray*}
where $\beta=\lambda(1+\langle \G(\sigma),\sigma \rangle)+\langle \sigma,\omega \rangle^2$. 
\end{corollary}




Moreover, if we consider  $\sigma_i\in\CV$, $i=1,\ldots,m+\ell$ such that $\sigma_i\notin \omega^\bot$ for  $i=1,\dots,m$ and $\sigma_{m+i}\in \omega^\bot$, $i=1,\ldots,\ell$  the operator
$$\H=\F+\sum_{i=1}^{m+\ell}\P_{\sigma_i},$$
is a perturbed operator. Here, $m$ or $\ell $ can be equal to $0$. The relation between the corresponding inverse operators  is given in the following theorem.

\begin{theorem}\cite[Theorem 3.5]{CEM14} \label{theorem35}  The operator $\H$ is positive semi-definite, positive definite when $m\ge 1$ and moreover
$${\mathcal H}^{\dag}=
\displaystyle \G+ h \P_{\omega}+\sum_{i=1}^{m+\ell} h_i [\P_{\G(\sigma_i),\omega}-
\P_{\omega,\G(\sigma_i)}]-\sum_{i,j=1}^{m+\ell} h_{ij}\P_{\G(\sigma_i),\G(\sigma_j)},$$
where   $(b_{ij})=(I+\langle\G(\sigma_j),\sigma_i\rangle)^{-1}$, and
$$\begin{array}{rl}
h=& \hspace{-.25cm}\displaystyle \left(\lambda+\sum_{r,s=1}^{m} b_{rs}\langle\sigma_r,\omega\rangle\langle\sigma_s,\omega\rangle\right)^{\dag}, \\[2ex]
h_i =& \hspace{-.25cm}h \displaystyle\sum_{r=1}^{m} b_{ir}\langle\sigma_r,\omega\rangle, \hspace{.25cm}i=1,\ldots,m+\ell,\\[2ex]
h_{ij}=&\hspace{-.25cm}\displaystyle b_{ij}-h\left(\sum_{r=1}^{m}b_{ir}\langle\sigma_r,\omega\rangle \right)\left(\sum_{r=1}^{m} b_{jr}\langle\sigma_r,\omega\rangle \right), \hspace{.25cm}i,j=1,\ldots,m+\ell.
\end{array}$$
\end{theorem}



On the other hand, the {\it Schur complement} provides us with a fundamental tool for the results on next sections.

\begin{lemma} \label{lemma1}
If $\AA\in \mathcal{M}_{n\times n}$, $\BB\in \mathcal{M}_{n\times m}$, $\DD\in \mathcal{M}_{m\times m}$ invertible, and $\Ss=\AA-\BB\DD^{-1}\BB^{\intercal}$, then
$$
 \left(\begin{array}{cc}
\AA & \BB\\
\BB^{\intercal} & \DD
\end{array}
\right)^{\dag}
=
\left(
\begin{array}{cc}
 \Ss^{\dag} & -\Ss^{\dag}\BB\DD^{-1}\\
-\DD^{-1}\BB^{\top}\Ss^{\dag}& \DD^{-1}+\DD^{-1}\BB^{\intercal}\Ss^{\dag}\BB\DD^{-1}
\end{array}
\right),$$
where $\displaystyle C^{\dag}$ stands for the Moore-Penrose inverse of the matrix $C$.

\end{lemma}

\section{Adding a new vertex}
\label{}


In the following, the triple $\Gamma=(V,E,c)$ denotes a finite network, i.e., a connected graph without loops nor multiple edges with vertex set $V$, with cardinality $n$, and edge set $E$, in which each edge $e_{xy}\in E$ has assigned a value $c(x,y)>0$ named conductance. The conductance $c$ is a symmetric function $c: V\times V \rightarrow [0,\infty)$ such that $c(x,x)=0$ for any $x\in V$ and where the vertex $x$ is adjacent to vertex $y$ iff $c(x,y)>0$.
Given a weight on $V$, $\omega\in \Omega(V)$, for any pair of vertices $(x,y)\in V\times V$ the $\omega$-dipole between $x$ and $y$ is the function $\tau_{xy}=\frac{\varepsilon_x}{\omega(x)}-\frac{\varepsilon_y}{\omega(y)}$. The  Laplacian of the network $\Gamma$
is the endomorphism of $\CV$ that assigns to each $u\in \CV$ the function
$$\L(u)(x)=\sum_{y\in V}c(x,y)[u(x)-u(y)], \quad x\in V.$$
The Laplacian is a singular elliptic operator on $\CV$ and moreover $\L(u)=0$ iff $u$ is a constant function. Given $q\in \CV$ the Schr\"odinger operator on $\Gamma$ with potential $q$ is the endomorphism of $\CV$ that assigns to each $u\in \CV$ the function $\L_q(u)=\L(u)+qu$. Given a weight $\omega\in \Omega(V)$, the potential determined by $\omega$ is the function $q_{\omega}=-\frac{1}{\omega}\L(\omega)$. It is well-known that the Schr\"odinger operator $\L_q$ is $(\lambda,\omega)$-elliptic iff $q=q_{\omega}+\lambda$, see \cite{BCE12}. Moreover, it is  singular iff $\lambda=0$ and then, $\L_{q_\omega}(v)=0$ iff $v=a\omega$, $a\in \mathbb{R}$. We denote by $\go$ the orthogonal Green operator associated with $\L_q$ and by $\sf{G}_{\lambda,\omega}$ its corresponding kernel. From now on, we consider fixed the value $\lambda\geq 0$, the weight $\omega \in \Omega(V)$ and the Schr\"odinger operator $\Lq$ with $q=q_{\omega}+\lambda$.




In this study we worry about perturbations of $\L_q$ performed by adding a new vertex. Namely, let $x'$ be a new vertex and assume we connect it to $m$ vertices $x_1,\dots, x_m\in V$, where $1\leq m \leq n$. The new network $\Gamma'=(V',E',c')$ has vertex set $V'=V\cup \{x'\}$, edge set $E'=E\cup \{e_{x_1x'},\dots,e_{x_mx'}\}$ and  conductance
$$c'(x,y)=
\left\{\begin{array}{ll}
a_i>0 & \textrm{if} \; \;x=x_i\in V,\; y=x', \quad 1\leq i \leq m,\\
c(x,y) & \textrm{if}  \; \;x,y\in V,\\
0& \textrm{otherwise}.
\end{array}
\right.
 $$

If $\omega(x')$ is a positive value assigned to $x'$, we define a  weight on $\Gamma'$ by $\omega'(x)=\omega(x)/\sqrt{1+\omega(x')^2}$, for any $x\in V'$. Observe that in this case it holds $\omega(x)/\omega(y)=\omega'(x)/\omega'(y)$, for any $x$, $y\in V'$.

The following notation is useful in what follows. For any $i=1,\ldots,m$, we denote by $\rho_i=\sqrt{a_i\omega(x_i)\omega(x')}$, $\sigma_i=\dfrac{\rho_i}{\omega(x_i)}\varepsilon_{x_i}$, $\sigma=\sum_{i=1}^m a_{i}\varepsilon_{x_i}$ and $\alpha=\lambda +\dfrac{1}{\omega^2(x')}\sum_{i=1}^m \rho_i^2$.

\begin{prop}
If $\L'$ is the Laplacian of $\Gamma'$ and $p=q_{\omega'}+\lambda,$  where $q_{\omega'}=-\dfrac{\L'(\omega')}{\omega'},$ then
$$\begin{array}{rll}
\L_p'=&\hspace{-.25cm}\L_q+\sum\limits_{i=1}^m\P_{\sigma_i}-\P_{\sigma,\varepsilon_{x'}} &\mbox{ on } V,\\[1ex]
\L_p'=&\hspace{-.25cm}-\P_{\varepsilon_{x'},\sigma}+\alpha\P_{\varepsilon_{x'}} &\mbox{ on } \{x'\}.
\end{array}$$
\end{prop}
\proof
For any $u\in\CVp$ it is satisfied that
$$\begin{array}{rl}
\L'u(x_i)=&\hspace{-.25cm}\L u(x_i)+a_i(u(x_i)-u(x')), \hspace{.25cm} i=1,\ldots,m,\\[1ex]
\L'u(x')=&\hspace{-.25cm}\sum\limits_{i=1}^ma_i(u(x')-u(x_i)),\\[1ex]
\L'u(x)=&\hspace{-.25cm}\L u(x), \hspace{.25cm} \hbox{ otherwise,}
\end{array}$$
and in particular,
$$\begin{array}{rl}
q_{\omega'}(x_i)=&\hspace{-.25cm}q_\omega(x_i)-a_i+a_i\dfrac{\omega(x')}{\omega(x_i)}, \hspace{.25cm} i=1,\ldots,m,\\[1ex]
q_{\omega'}(x')=&\hspace{-.25cm}-\sum\limits_{i=1}^ma_i +\dfrac{1}{\omega(x')}\sum\limits_{i=1}^ma_i\omega(x_i),\\[1ex]
q_{\omega'}(x)=&\hspace{-.25cm}q_\omega(x), \hspace{.25cm} \hbox{ otherwise}.
\end{array}$$
Therefore,
$$\begin{array}{rl}
\L'_pu(x_i)=&\hspace{-.25cm}\L_q u(x_i)+\dfrac{\rho_i^2}{\omega^2(x_i)}u(x_i)-a_iu(x'), \hspace{.25cm} i=1,\ldots,m,\\[3ex]
\L'_pu(x')=&\hspace{-.25cm}\displaystyle \lambda u(x')+\dfrac{u(x')}{\omega^2(x')}\sum\limits_{i=1}^m\rho_i^2-\sum\limits_{i=1}^ma_iu(x_i),\\[3ex]
\L'_pu(x)=&\hspace{-.25cm}\L_q u(x), \hspace{.25cm} \hbox{ otherwise},
\end{array}$$
and the result follows.\qed

 The relation between the matrices associated with both Schr\"odinger operators of $\Gamma$ and $\Gamma'$ is given by
$$ {\sf L}'_{p}=
\left(\begin{array}{cc}
{\sf H} & -{\sf s}\\
 -{\sf s}^{\intercal} & \alpha
\end{array}
\right),$$
where ${\sf H}$ is the matrix associated with  the  operator
$$\mathcal{H}=\Lq+\sum_{i=1}^m\P_{\sigma_i},$$
and ${\sf s}$ is the column vector ${\sf s}=\sum_{i=1}^m a_{i}{\sf e}_{x_i}$, where ${\sf e}_{i}$ for $i=1,\dots,n$ are the vectors of the canonical basis.

In order to compute the Moore--Penrose inverse of ${\sf L}'_{p}$ we will use Lemma \ref{lemma1} and it will be useful to introduce the following perturbations. We define, for $k=1,\dots,m,$
$$
 \pi_k=\sqrt{\dfrac{\lambda}{\alpha}}\sigma_k,$$
 and for $i=1,\ldots,m-1$, $j=i+1,\ldots,m,$ let $k=\dfrac{(2m-1-i)i}{2}+j$ and
 $$ \pi_k=\dfrac{1}{\sqrt{\alpha}}\dfrac{\rho_i\rho_j}{\omega(x')}\left(\dfrac{\varepsilon_{x_i}}{\omega(x_i)}-\dfrac{\varepsilon_{x_j}}{\omega(x_j)}\right).$$

\begin{theorem} \label{maintheorem}
 The Moore--Penrose inverse of ${\sf L}'_{p}$ is given by
$$
 ({\sf L}'_{p})^{\dag}
=
\left(
\begin{array}{cc}
 {\sf M} & -\frac{1}{\alpha}{\sf M}{\sf s}\\[1ex]
-\frac{1}{\alpha}{\sf s}^{\intercal}{\sf M}& \frac{1}{\alpha}+\frac{1}{\alpha^2}{\sf s}^{\intercal}{\sf M}{\sf s}
\end{array}
\right),
$$
where ${\sf M}$ is the matrix associated with the operator
$$\displaystyle \G+ h \P_{\omega}+\sum_{i=1}^{\frac{m(m+1)}{2}} h_i [\P_{\G(\pi_i),\omega}-
\P_{\omega,\G(\pi_i)}]-\sum_{i,j=1}^{\frac{m(m+1)}{2}} h_{ij}\P_{\G(\pi_i),\G(\pi_j)},$$
and  where if  $(b_{ij})=(I+\langle\G(\pi_j),\pi_i\rangle)^{-1}$,
$$\begin{array}{rl}
h=& \hspace{-.25cm}\displaystyle \lambda^\dag\alpha \left(\alpha+\sum_{r,s=1}^{m} b_{rs}\rho_r\rho_s\right)^{-1}, \\[2ex]
h_i =& \hspace{-.25cm}h \sqrt{\dfrac{\lambda}{\alpha}}\displaystyle\sum_{r=1}^{m} b_{ir}\rho_r, \hspace{.25cm}i=1,\ldots,\frac{m(m+1)}{2},\\[2ex]
h_{ij}=&\hspace{-.25cm}\displaystyle b_{ij}-\dfrac{h\lambda}{\alpha}\left(\sum_{r=1}^{m}b_{ir}\rho_r \right)\left(\sum_{s=1}^{m} b_{js}\rho_s\right), \hspace{.25cm}i,j=1,\ldots,\frac{m(m+1)}{2}.
\end{array}$$
\end{theorem}
\proof
From Lemma \ref{lemma1}, we have that
$$
({\sf L}'_{p})^{\dag}=\left(
\begin{array}{cc}
 \Ss^{\dag} & -\dfrac{1}{\alpha}\Ss^{\dag}{\sf s}\\
-\dfrac{1}{\alpha}{\sf s}^\intercal\Ss^{\dag}& \frac{1}{\alpha}+\frac{1}{\alpha^2}{\sf s}^{\intercal}{\sf S}^{\dag}{\sf s}
\end{array}
\right),
$$
where
$\Ss={\sf H}-\dfrac{1}{\alpha}{\sf s}\otimes{\sf s}$ is the matrix associated with the operator %
$${\mathcal S}=\Lq+\sum_{i=1}^m\P_{\sigma_i}-\dfrac{1}{\alpha}\P_\sigma.$$
Now, let us prove that
$$\P_{\sigma}=(\alpha-\lambda)\sum_{i=1}^m\P_{\sigma_i}-\sum_{1\le i<j\le m}\P_{\sigma_{ij}},$$
where $\sigma_{ij}=\sqrt{\alpha}\pi_k,$  $k=\dfrac{(2m-1-i)i}{2}+j$ for $i=1,\ldots,m-1$, $j=i+1,\ldots,m.$  If  $P_\sigma$ denotes the kernel of $\P_\sigma$ and $P$ denotes the kernel of the operator in the right side of the equality,  the claim is equivalent  to prove that $P=P_\sigma$.

Since $P_\sigma=\sigma\otimes \sigma$, we have
$$P_\sigma=\sum\limits_{i,j=1}^ma_ia_j(\varepsilon_{x_i}\otimes \varepsilon_{x_j}).$$

On the other hand, for any $i=1,\ldots,m$, we have
$$P_{\sigma_i}=\dfrac{\rho_i^2}{\omega(x_i)^2}\,(\varepsilon_{x_i}\otimes \varepsilon_{x_i})=\dfrac{a_i\omega(x')}{\omega(x_i)}\,(\varepsilon_{x_i}\otimes \varepsilon_{x_i}),$$
and hence,
$$\sum\limits_{i=1}^mP_{\sigma_i}=\omega(x')\sum\limits_{i=1}^m\dfrac{a_i}{\omega(x_i)}\,(\varepsilon_{x_i}\otimes \varepsilon_{x_i}).$$

Moreover, for $1\le i<j\le m$,
$$\begin{array}{rl}
P_{\sigma_{ij}}=&\hspace{-.25cm}\displaystyle \dfrac{\rho_i^2\rho_j^2}{\omega(x')^2\omega(x_i)^2}\,(\varepsilon_{x_i}\otimes \varepsilon_{x_i})+\dfrac{\rho_i^2\rho_j^2}{\omega(x')^2\omega(x_j)^2}\,(\varepsilon_{x_j}\otimes \varepsilon_{x_j})\\[3ex]
-&\hspace{-.25cm}\displaystyle\dfrac{\rho_i^2\rho_j^2}{\omega(x')^2\omega(x_i)\omega(x_j)}\big(\varepsilon_{x_i}\otimes \varepsilon_{x_j}+\varepsilon_{x_j}\otimes \varepsilon_{x_i}\big)\\[3ex]
=&\hspace{-.25cm}\displaystyle \dfrac{a_ia_j\omega(x_j)}{\omega(x_i)}\,(\varepsilon_{x_i}\otimes \varepsilon_{x_i})+\dfrac{a_ia_j\omega(x_i)}{\omega(x_j)}\,(\varepsilon_{x_j}\otimes \varepsilon_{x_j})\\[3ex]
-&\hspace{-.25cm}\displaystyle a_ia_j \big(\varepsilon_{x_i}\otimes \varepsilon_{x_j}+\varepsilon_{x_j}\otimes \varepsilon_{x_i}\big),
\end{array}$$
and hence,
$$\begin{array}{rl}
\displaystyle \sum\limits_{1\le i<j\le m}P_{\sigma_{ij}}
=&\hspace{-.25cm}\displaystyle \dfrac{1}{2} \sum\limits_{i=1}^m\sum\limits_{j=1\atop j\not=i}^m\dfrac{a_ia_j\omega(x_j)}{\omega(x_i)}\,(\varepsilon_{x_i}\otimes \varepsilon_{x_i})+\dfrac{1}{2} \sum\limits_{j=1}^m\sum\limits_{i=1\atop i\not=j}^m\dfrac{a_ia_j\omega(x_i)}{\omega(x_j)}\,(\varepsilon_{x_j}\otimes \varepsilon_{x_j})\\[3ex]
-&\hspace{-.25cm}\displaystyle\displaystyle \sum\limits_{1\le i<j\le m} a_ia_j \big(\varepsilon_{x_i}\otimes \varepsilon_{x_j}+\varepsilon_{x_j}\otimes \varepsilon_{x_i}\big)\\[3ex]
=&\hspace{-.25cm}\displaystyle  \sum\limits_{i=1}^m\dfrac{a_i}{\omega(x_i)}\Big(\sum\limits_{j=1\atop j\not=i}^ma_j\omega(x_j)\Big)(\varepsilon_{x_i}\otimes \varepsilon_{x_i})- \sum\limits_{i,j=1\atop i\not=j}^ma_ia_j\big(\varepsilon_{x_i}\otimes \varepsilon_{x_j}\big).\end{array}$$

Taking into account that
$$\sum\limits_{j=1\atop j\not=i}^ma_j\omega(x_j)=\sum\limits_{j=1}^ma_j\omega(x_j)-a_i\omega(x_i)=(\alpha-\lambda)\omega(x')-a_i\omega(x_i),$$
we obtain that
$$\sum\limits_{1\le i<j\le m}P_{\sigma_{ij}}
=(\alpha-\lambda)\sum\limits_{i=1}^mP_{\sigma_i}-\sum\limits_{i,j=1}^ma_ia_j\big(\varepsilon_{x_i}\otimes \varepsilon_{x_j}\big).
$$
Therefore,
$$
{\mathcal S}=\Lq+\sum\limits_{k=1}^{\frac{m(m+1)}{2}}\P_{\pi_k}.$$
Finally, from Theorem \ref{theorem35}, we get that ${\Ss}^\dag={\sf M}$ and the result follows.\qed

Next we describe the coefficients of matrix $(\langle \G(\pi_\ell),\pi_k \rangle)$.

\begin{lemma}         \label{matrixperturbations}
The elements of the matrix $A=(\langle \G(\pi_\ell),\pi_k \rangle)$ are given by
$$\begin{array}{rl}
(A)_{k,	\ell}=&\hspace{-.25cm}\displaystyle \frac{\lambda\rho_k\rho_\ell}{\alpha}{\dfrac{{\sf G}_{\lambda,\omega}(x_k,x_\ell)}{\omega(x_k)\omega(x_\ell)}}, \hspace{.25cm} k,\ell=1,\ldots,m,\\[3ex]
(A)_{k,\ell}=&\hspace{-.25cm}\displaystyle \frac{\sqrt{\lambda}\rho_k\rho_i\rho_j}{\alpha\omega(x')}\left[\frac{{\sf G_{\lambda,\omega}}(x_k,x_i)}{\omega(x_k)\omega(x_i)}-
\frac{{\sf G_{\lambda,\omega}}(x_k,x_j)}{\omega(x_k)\omega(x_j)}\right],\\[3ex]
            & k=1,\dots,m,\; \ell=\dfrac{(2m-1-i)i}{2}+j,\\[3ex]
(A)_{k,\ell}=&\hspace{-.25cm}\displaystyle \frac{\rho_i\rho_j\rho_r\rho_s}{\alpha\omega(x')^2}
        \left[\frac{{\sf G_{\lambda,\omega}}(x_i,x_r)}{\omega(x_i)\omega(x_r)}-\frac{{\sf G_{\lambda,\omega}}(x_i,x_s)}{\omega(x_i)\omega(x_s)}-\frac{{\sf G_{\lambda,\omega}}(x_j,x_r)}{\omega(x_j)\omega(x_r)}+\frac{{\sf G_{\lambda,\omega}}(x_j,x_s)}{\omega(x_j)\omega(x_s)}\right],\\[3ex]
         & \hspace{.25cm}  for \; r=1,\ldots,m-1,\;s=r+1,\ldots,m,\;  k=\dfrac{(2m-1-r)r}{2}+s, \\[3ex]
         & \hspace{.25cm}  and \;for\; i=1,\ldots,m-1,\; j=i+1,\ldots,m,\;  \ell=\dfrac{(2m-1-i)i}{2}+j.
\end{array}$$
\end{lemma}

Observe that we can deduce two special cases of  Theorem \ref{maintheorem}; if $m=1$ that means the addition of a pendant vertex to the network (see~\cite{CEGM14}), and if $m=n$ it represents the join of a the new vertex with the graph; (see~\cite{BCE12}).




\begin{corollary} \label{corolpendant}If $m=1$ and $a=c(x,x')$, then
$$
 ({\sf L}'_{p})^{\dag}
=
\left(
\begin{array}{cc}
 {\sf M} & -\frac{a}{\alpha}{\sf M}{\sf e}_x \\[1ex]
-\frac{a}{\alpha}{\sf e}_x ^\intercal{\sf M}& \frac{1}{\alpha}+\frac{a^2}{\alpha^2} {\sf e}_x ^\intercal{\sf M}{\sf e}_x
\end{array}
\right),
$$
where ${\sf M}$ is the matrix associated with the operator
$$
\G- \frac{1}{h}
\Big[\lambda \P_{\G(\sigma)} -\rho_x \left(\P_{\G(\sigma),\omega} -  \P_{\omega,\G(\sigma)}\right)  -\left(1+ (\alpha-\lambda){\sf G}_{\lambda,\omega}(x,x)\right)\P_{\omega} \Big],
$$
where $h=\lambda[1+(\alpha-\lambda){\sf G}_{\lambda,\omega}(x,x)]+\rho_x^2$.

\end{corollary}

%


\section{Acknowledgement}
This work has been supported by the Spanish Research Council (Ministerio de Ciencia e Innovaci\'on) under the projects MTM2011-28800-C02-02 and MTM2011-28800-C02-01.



\bibliographystyle{model1-num-names}
\bibliography{<your-bib-database>}



\end{document}